\documentclass{amsart}
\usepackage{amssymb}

\newcommand{\beq}[1]{\begin{equation}\label{#1}}
\newcommand{\eeq}{\end{equation}}

\theoremstyle{definition}

\title{ Ramsey-product subsets of a group }

\author{  Igor Protasov and Ksenia Protasova }

\begin{document}


\begin{abstract} 
We say that a subset $S$ of an infinite group $G$ is a Ramsey-product subset if, for any infinite subsets $X$, $Y$  of $G$, there exist $x \in X$  and $y\in Y$ such that $x y \in S$  and  $ y x \in S$ .  We show that the family $\varphi$ of all  Ramsey-product subsets of $G$ is a filter and $\varphi$ defines the subsemigroup $ \overline{G^*G^*}$ of the semigroup $G^*$ of all free ultrafilters on $G$. 
\end{abstract}
\maketitle

\Large

All groups under consideration are supposed to be infinite; a countable set means a    countably infinite set.

We say that a subset $S$ of a group $G$ is  \vskip 5pt

\begin{itemize}

\item{}   a {\it Ramsey-product subset }  if, for any infinite subsets $X, Y$  of $G$, there exist  $x\in X$  and $y\in Y$  such that $xy\in S $ and  $yx\in S $; \vskip 5pt

\item{}   a {\it Ramsey-square subset }  if, for every infinite subset $X$  of $G$ there exist distinct elements $x, y\in X$   such that $xy\in S $,   $yx\in S $; \vskip 5pt

\item{}    a  {\it Ramsey-quotient  subset }    if,  for every infinite subset $X$ of $G$, there exist  distinct elements  $x,y\in X$  such  that $x^{-1} y\in S $  and $ y^{-1} x \in S $.

\end{itemize}
\vskip 10pt
\Large

We show that above defined subsets of $G$ arise naturally in context of the   Stone-$\check{C}$ech compactification  $\beta G$  of the group $G$ endowed with the discrete topology. 
We  identify $\beta G$  with the set of all ultrafilters on $G$.  
Then the family $\{ \bar{A}: A\subseteq  G \} $, where $ \bar{A}=\{ p \in  \beta G: A\in p\}$,   forms a base for the topology of $\beta G$.  Given a filter $\varphi $  on $G$, we denote  $ \bar{\varphi} = \cap \{  \bar{A} : A\in \varphi \}$, so  $ \varphi$ defines the closed subset $ \bar{\varphi}$  of $\beta G$,
and every non-empty closed subset of $\beta G$ can be defined in this way . 

We use the standard extension [3, Section 4.1] of the multiplication on $G$ to the semigroup multiplication on $\beta G$. Given two ultrafilters $p,q \in \beta G$,  we choose $P\in p$  and, for each $x\in P$,  pick $Q _x \in q$. 
Then $\bigcup _{x\in P}  xQ_x \in pq$  and the family of these subsets forms a base of the product $pq$.  We note that the set $G^*$  of all free ultrafilters of $G$ is a closed subsemigroup of $\beta G$, and the closure $ \overline{G^*G^*}$  of $G^*G^*$ is a subsemigroup of $G^*$. An ultrafilter $p\in G^*\setminus \overline{G^*G^*}$ is called {\it strongly prime}.

By [5, Propositions 3, 4], the family $\psi $ of all  Ramsey-quotient  subsets   
of a group $G$  is a filter and $\overline{\psi}$ is the smallest closed subset of $G ^*  $  containing all ultrafilters of the form $q^{-1} q$, $q\in G^*$, where $q^{-1}=\{ Q^{-1} : Q\in q \}$. Analogously, the closure of $\{pp: p\in G^* \}$   is defined by the filter of all Ramsey-square subsets of $G$. In both cases, we used the classical Ramsey   theorem [2, p. 16].  

In this note, the key technical part plays the following  version of 
Ramsey  theorem (Theorem 6 from [2, p.98]).
\vskip 5pt

{\it Let $\chi$ be a finite coloring $\chi : \omega  \times   \omega  \rightarrow [r]$. Then there exists an infinite set  $A=\{ a_{i} \}_{< \omega}$ and colors $c_L$,   $c_G$,  $c_E$ (not necessarily distinct) such that $\chi (a_i, a_j)=c_L$ if $i<j$ ,  $\chi (a_i, a_j)=c_G$   if  $i>j$  and    $\chi (a_i, a_j)=c_E$ if  $i=j$.}
\vskip 15pt

{\bf  Lemma 1}. {\it If $S$ is a Ramsey-product subsets of a group $G$ then, for any countable subsets  $X, Y$ of $G$, there exist disjoint countable subsets $X '  \subseteq  X$  and $Y '  \subseteq  Y$ such that  $X ' Y ' \subseteq  S$,   $ Y ' X ' \subseteq  S$.
\vskip 10pt

Proof.} Passing to subsets of $X$ and $Y$, we may suppose that $X\cap Y=\emptyset$. We enumerate $X=\{ x_n : n< \omega\}$, $Y=\{ y_n : n<\omega \}$ and define a $\{0,1\}$-coloring $\chi$   of  $\omega \times \omega$ by the rule: $\chi ((n,m))= \chi ((m,n))=1$ if $m>n$,  $x_n y_m \in S $, $ y_m x_n \in S$, and $\chi ((m,n))=0$  otherwise. By Theorem 6 from [2, p.98], there exists a countable subset $A$ of $\omega$ such that $\chi$ is monochrome on all $(n,m)$ , such that $n\in A$, $m\in A$, $n\neq   m$. We partition $A$  into two infinite subsets $A=B\cup C$ and denote $X ' = \{ x_n : n\in B\}$,  $Y ' = \{ y_n : n\in C\}$. Since $S$ is a Ramsey-product set, there are $x' \in X ' $, 
$y ' \in Y ' $ such that $x'  y' \in S$, $y ' x ' \in S $. By the choice of $A$, we have $xy\in S$, $yx\in S$ for all $x\in X '$, $y\in Y ' $.   $ \ \  \ \ \Box$

\vskip 15pt

{\bf  Lemma 2}. {\it For  a free ultrafilter $r$ on a group 
$G$,  $r\in \overline{G^* G^*}$ if and only if, for every subset $R\in r$, there exist two countable subsets $\{x_n : n<\omega  \}$  and $\{y_n : n<\omega  \}$   of $G$ such that $\{x_n  y_m : n  \leq     m <\omega  \} \subseteq R.$  
\vskip 10pt

Proof.} We assume that $r\in  \overline{G^* G^*}$ , take an arbitrary $R\in r$  and choose $p,q \in G*$ such that $R\in pq$. We choose $P\in p$  and $Q_x$, $x\in P$  such that $\bigcup _{x\in P} x Q_x \subseteq R$. We take an arbitrary countable subset $\{ x_n : n<\omega  \}$ of $P$ and choose an injective sequence $(y_n)_{n\in \omega}$ in $G$  such that $y_n \in Q_{x_0} \cap\ldots \cap Q_{x_n} $, $n<\omega$. Then $\{ x_n y_m :  n\leq   m <\omega \}\subseteq R$.

On the other hand, let $R\in r$  and $\{ x_n : n\in \omega \}$,  $\{ y_n : n\in \omega \}$ are chosen so that $\{ x_n  y_m : n  \leq  m < \omega \} \subseteq R$.  We take an arbitrary free ultrafilters $p,q$ such that $\{ x_n : n\in \omega \}\in P$,  $\{ y_n : n\in \omega \}\in q$.  Then $\{ x_n  y_m : n  \leq  m < \omega \}\in pq$, so $R\in pq$  and $r\in \overline{G^* G^*}$.     $ \ \  \ \ \Box$

\vskip 15pt

We recall [1]  that a subset $A$  of a group $G$ is {\it sparse} if, for  every infinite subset $X$ of $G$, there exists a non-empty finite subset $F\subset X$ such that $\bigcap _{x\in F} xA$ is finite. By [1, Theorem 9], an ultrafilter $r\in G^*$ is strongly prime if  and only if there exists a sparse subset $A$ of $G$ such that $A\in r$. For sparse subsets see also [4], [6]. 

\vskip 15pt 

{\bf  Teorem}. {\it  For every group $G$, the following statements hold: \vskip 5pt

$(i)$  the family $\varphi$  of all Ramsey-product subsets of  $G$ is a filter;\vskip 5pt

$(ii)$  $ \overline{\varphi } =  \overline{G^* G^*} $;\vskip 5pt

$(iii)$  $A\in \varphi$  if and only if $G\setminus A$ is sparse.

\vskip 10pt

Proof.} $(i)$ We take two Ramsey-product subsets $S$ and $T$ of $G$  and show that $S\cap T\in \varphi .$ 
Let $X,Y$ be infinite subsets of $G$. 
By Lemma 1, there exist disjoint countable subsets $X' \subseteq X$  and $Y' \subseteq Y$ such that $xy\in S$,  $yx\in S$ for all  $x\in X'$,  $y\in Y'$. Since $T$ is a Ramsey-product subsets, there exist $x'\in X'$,  $y'\in Y'$
such that   $x'y'\in T$,  $y'x'\in T$.  Then $x'y'\in S \cap T$,  $y'x'\in S\cap T$, so $S\cap T\in  \varphi $.
\vskip 5pt

 $(ii)$  To show that $ \overline{G^* G^*} \subseteq   \overline{\varphi}$, we suppose the contrary, choose
$p\in  \overline{G^* G^*} \setminus   \overline{\varphi}$  and take $P\in p$ such that $G\setminus P\in \varphi $.  By Lemma 2, there are countable sets $\{ x_n : n<\omega \}$ and  $\{ y_n : n<\omega \}$ such that  $\{ x_n y_m : n\leq  m<\omega \}\subseteq P $.  We use Lemma 1 to choose disjoint countable subsets  $X' \subseteq \{ x_n : n<\omega \} $,  $Y' \subseteq \{y_n : n<\omega \} $ such that 
 $X' Y' \subseteq  G\setminus P$,   $ \  \  Y'X' \subseteq  G\setminus P$. To get  a contradiction, we take an arbitrary $x_n \in X'$  and pick $y_m \in Y'$ such that $m>n$. 

Now we take an arbitrary   $p\in   \overline{\varphi} $  and prove that  $p\in \overline{G^* G^*}$. Given any $P\in p$, by Lemma 1, there exist two disjoint countable subsets $X, Y$  of $G$ such that, for any countable subsets $X'\subseteq  X$  and  $Y'\subseteq  X$,   either   $X'Y'\cap P \neq  \emptyset $
  or $Y'X'\cap P \neq  \emptyset $ . 
We enumerate  $X=\{ x_n: n<\omega  \}$,   $Y=\{ y_n: n<\omega  \}$  and define a $\{0,1\}$-coloring $\chi$  of   $\omega\times\omega$  by the rule:  $\chi(n,m)=1$
 if and only if  $n<m$ and $x_n y_m\in P$ or  $n>m$ and  $y_m x_n \in P$.
 Applying Theorem 6 from [2, p. 98], we get an infinite subset $W$ of $\omega$ such that the restriction of $\chi$ to $W_1=\{ (n, m):  n<m ,  \   \  n\in W,  \   \  m\in W \} $
  is monochrome, and the restriction of $\chi$ to $W_2=\{ (n, m):  n>m ,  \   \  n\in W,  \   \  m\in W \} $
 is monochrome. 
By the choice of $X$  and  $Y$, either  $\chi |_{W_{1}} \equiv  1$
  or    $\chi |_{W_{2}} \equiv  1$.
In the first case, we choose two injective sequences  $(x_n ')_{n<\omega} $   in $X$ and  $(y_n ')_{n<\omega} $   in $Y$  such that  $\{ x_{n} '  y_{m} ' : n<m< \omega \}\subseteq P  $ . 
In the second case, we choose two injective sequences
$(x_n ')_{n<\omega} $ 
 in $X$  and  $(y_n ')_{n<\omega} $
 in $Y$ such that 
$\{ y_n ' x_m '  : n<m< \omega \}\subseteq P  $ . 
Applying Lemma 2, we conclude that   $p\in  \overline{G^* G^*}   $. 

\vskip 5pt
$(iii)$   We apply $(ii)$ and Theorem 9  from [1].   $ \ \  \ \ \Box$

\vskip 15pt

We recall that a subset $S$ of a group $G$ is {\it  large (extralarge) } if there is a finite subset $K$  of $G$  such that $G=KS$  ($S\cap L$  is large for every  large  subset   $L$ of $G$).  By Theorem 4.1 from [4] , the complement of sparse subset is extralarge. Applying $(iii)$, we see that each   Ramsey-product subset is extralarge.

By the statement  $(ii)$, every Ramsey-product subset of $G$ is a member of every idempotent of the semigroup $G^*$, in particular, every  idempotent from  the minimal ideal of  $\beta G$,  so $S$ is combinatorially  very rich (see [3, Section 14]).

If $G$ is an amenable group and $A$ is a sparse subset of $G$, by Theorem 5.1 from [4], $\mu (A)=0$  for every left invariant  Banach measure  $\mu$  on $G$. Applying (iii), we see that $\mu (S)=1$
for every Ramsey-product subset $S$ of $G$ and every left invariant Banach measure $\mu$ on $G$.

We conclude the note with the following formal generalization. We fix two infinite subset $X, Y$  of a group $G$ and say that a subset $S$ of $G$ is a 
{\it Ramsey $(X, Y)$-product subset } if, for any infinite subsets   $X'\subseteq  X$ and $Y'\subseteq  Y$, there exist  $x\in X'$  and $y\in Y'$  such that $xy\in S$, $yx\in S$. Then the family $\varphi _{X,Y}$ of all Ramsey 
 $(X, Y)$-product subsets of $G$ is a filter and $ \overline{\varphi} _{X,Y}= \overline{X^* Y^*}  \   \cup \   \overline{Y^* X^*} $.

\vspace{5 mm}
CONTACT INFORMATION

I.~Protasov: \\
Faculty of Computer Science and Cybernetics  \\
        Kyiv University  \\
         Academic Glushkov pr. 4d  \\
         03680 Kiev, Ukraine \\ i.v.protasov@gmail.com

\medskip

K.~Protasova:\\
Faculty of Computer Science and Cybernetics \\
        Kyiv University  \\
         Academic Glushkov pr. 4d  \\
         03680 Kiev, Ukraine \\ ksuha@freenet.com.ua

\end{document}